%% file: Coefficients.tex
\documentclass[review,onefignum,onetabnum]{siamart171218}

\input{ex_shared}
\nolinenumbers

\ifpdf
\hypersetup{
  pdftitle={On the coefficients in an asymptotic expansion of $(1+1/x)^x$},
  pdfauthor={T. M. Dunster and Jessica M. Perez}
}
\fi

\begin{document}

\maketitle

\begin{abstract}
 The function $g(x)= (1+1/x)^{x}$ has the well-known limit $e$ as $x\rightarrow{\infty}$. The coefficients $c_{j}$ in an asymptotic expansion for $g(x)$ are considered. A simple recursion formula is derived, and then using Cauchy's integral formula the coefficients are approximated for large $j$. From this it is shown that $|c_{j}|\rightarrow{1}$ as $j\rightarrow{\infty}$.
\end{abstract}

\begin{keywords}
  {Series expansions, Asymptotic representations in the complex plane, Integrals of Cauchy type}
\end{keywords}

\begin{AMS}
  41A58, 30E15, 30E20
\end{AMS}

\section{Introduction and main results}

Chen and Choi \cite{Chen:2014:AAF} presented a method to estimate Euler's constant $e$, accurate to as many decimal places as desired. Their starting point was the well known limit $\lim_{x\rightarrow \infty
}g(x)=e$, where
\begin{equation}
g(x)=\left(1+\frac{1}{x}\right) ^{x}.  
\label{g}
\end{equation}
Their method was based on an asymptotic expansion of $g(x)$ for large values of $x$, namely
\begin{equation}
g(x)\sim \sum_{j=0}^{\infty }\frac{c_{j}}{x^{j}}. \label{eu}
\end{equation}
We shall show that this asymptotic expansion
converges for $x>1$.

In \cite{Chen:2014:AAF} they proved that $c_{j}=e a_{j}$, where $a_{j}$ are rational numbers which alternate in sign, and are given explicitly by 
\begin{equation}
a_{j}=(-1)^{j}\ \sum \frac{1}{k_{1}!k_{2}!\cdots k_{j}!}\ \left( \frac{1}{2}
\right) ^{k_{1}}\left( \frac{1}{3}\right) ^{k_{2}}\cdots \left( \frac{1}{j+1}
\right) ^{k_{j}}.
\label{chen}
\end{equation}
Here for each $j$ the sum is taken over all possible combinations of nonnegative integers $k_{1},k_{2},k_{3},\cdots ,k_{j}$ that satisfy the
relation $\sum_{l=1}^{j} l k_{l}=j$. We remark that Ponomorenko \cite{Ponomarenko:2015:AFF} gave a much simpler proof of (\ref{chen}) using Fa\`{a} di Bruno's formula generalizing the chain rule to higher derivatives. 

The number of terms in the sum (\ref{chen}) is the partition function $P(j)$ \cite[Sect. 27.14(i)]{NIST:DLMF}. Hardy and Ramanujan \cite{Hardy:1918:AFI} gave the asymptotic formula 
\begin{equation}
P(j)\sim \frac{\exp (\pi \sqrt{2j/3})}{4j\sqrt{3}}
\end{equation}
as $j\rightarrow \infty $. It is therefore evident that the number of terms in the formula (\ref{chen}) grows exponentially in $\sqrt{j}$, and so is only practicable for small or moderate values of $j$. 

We derive a new way of computing the coefficients $c_{j}$, by a simple recursion formula. We also provide a simple asymptotic approximation for the coefficients as $j\rightarrow{\infty}$, and this shows that in absolute value they approach the value 1.

Our main results read as follows. In \cref{sec2} we prove (\ref{main2}), and in \cref{sec3} we prove (\ref{main3}).

\begin{theorem}
For $x>1$ the
expansion (\ref{eu}) converges, where the coefficients ${c_{j}}$ are given
recursively by $c_{0}=e$ and 
\begin{equation}
c_{j+1}=\frac{1}{j+1}{\sum_{l=0}^{j}(-1)^{j-l+1}\left( \frac{j-l+1}{j-l+2}
\right) c_{l}\quad (j=0,1,2,\cdots )}.  
\label{main2}
\end{equation}
Moreover as $j\rightarrow \infty $ 
\begin{equation}
c_{j}=(-1)^{j}\left(1+\frac{1}{j}\right) 
+\mathcal{O}\left(\frac{\ln(j)}{j^{2}}\right) .  
\label{main3}
\end{equation}
\end{theorem}
\begin{remark}
Since $\lim_{x\rightarrow \infty}g(x)=e$ it is clear from (\ref{eu}) that $c_{0}=e$. 
\end{remark}

\section{Proof of the recursion formula (\ref{main2})}
\label{sec2}
Define $z=1/x$ and $f(z)=g(1/z)$ so that (\ref{eu}) is written in the new form 
\begin{equation}
f(z)=(1+z)^{1/z}={\sum_{j=0}^{\infty }c_{j}z^{j}}.  
\label{1.8.1}
\end{equation}
and in this we consider $z$ complex. Moreover, on writing it as 
\begin{equation}
f(z)=\exp \left\{ z^{-1}\ln (1+z)\right\} ,  
\label{1.9.3}
\end{equation}
we note that it has a removable singularity at $z=0$, and therefore can be considered analytic at $z=0$ by assuming $f(0)=\lim_{z\rightarrow
0}f(z)=c_{0}=e$. We also see it has one finite singularity (a logarithmic branch point) at $z=-1$. Therefore the radius of convergence of the series (\ref{1.8.1}) is 1, i.e. it converges for $\left\vert z\right\vert <1$. Thus (\ref{eu}) converges for $x>1$, as asserted. By taking the principal branch of the logarithm in (\ref{1.9.3}) we have that $f(z)$ is analytic on the cut plane $\mathbb{C} \setminus (-\infty ,-1]$.

Next, on taking the natural logarithm of both sides of (\ref{1.8.1}), we get 
\begin{equation}
{\frac{1}{z}\ln (1+z)=\ln }\left\{ {\sum_{j=0}^{\infty }c_{j}z^{j}}\right\}.
\label{1.8.2}
\end{equation}
Then expand the $\ln (1+z)$ term by its Maclaurin expansion valid for $\left\vert z\right\vert <1$ and we arrive at 
\begin{equation}
\sum_{j=1}^{\infty }(-1)^{j}\frac{z^{j}}{j+1}=\ln \left\{ \sum_{j=0}^{\infty
}c_{j}z^{j}\right\} .  
\label{1.8.5}
\end{equation}
Next differentiate both sides with respect to $z$ to yield 
\begin{equation}
\sum_{j=1}^{\infty }(-1)^{j}\frac{jz^{j-1}}{j+1}=\sum_{j=1}^{\infty
}jc_{j}z^{j-1}\left[ \sum_{j=0}^{\infty }c_{j}z^{j-1}\right] ^{-1}.
\label{1.8.6}
\end{equation}
By shifting indices of the series starting at $j=1$ to start at $j=0$, and taking the series $\sum_{j=0}^{\infty }c_{j}z^{j}$ to the left-hand side, we see this is equivalent to 
\begin{equation}
\sum_{j=0}^{\infty }c_{j}z^{j}\sum_{j=0}^{\infty
}d_{j}z^{j}=\sum_{j=0}^{\infty }(j+1)c_{j+1}z^{j},  \label{1.8.7}
\end{equation}
where $d_{j}$ is given by 
\begin{equation}
d_{j}=(-1)^{j+1}\frac{j+1}{j+2}.  
\label{1.8.8}
\end{equation}

We now use the Cauchy product, which is the discrete convolution of two infinite series. It is given by the formula \cite[Sect. 73]{Brown:2014:CVA} 
\begin{equation}
\sum_{j=0}^{\infty }C_{j}\sum_{j=0}^{\infty }D_{j}=\sum_{j=0}^{\infty }\left[\sum_{l=0}^{j}C_{l}D_{j-l}\right].
\label{1.8.9}
\end{equation}
Applying this to the left-hand side of (\ref{1.8.7}) we combine both power series to the following single power series 
\begin{equation}
\sum_{j=0}^{\infty }c_{j}z^{j}\sum_{j=0}^{\infty
}d_{j}z^{j}=\sum_{j=0}^{\infty }\left[ {\sum_{l=0}^{j}c_{l}d_{j-l}}\right]
z^{j}{.}  \label{1.8.11}
\end{equation}
Finally substitute this into the left-hand side of (\ref{1.8.7}), and equate coefficients of $z^{j}$, to obtain 
\begin{equation}
{\sum_{l=0}^{j}c_{l}d_{j-l}}=(j+1)c_{j+1},  \label{1.8.13}
\end{equation}
and then using (\ref{1.8.8}) this leads to (\ref{main2}).

\section{Proof of the asymptotic approximation (\ref{main3})}
\label{sec3}
We will use the famous Cauchy integral formula  \cite[Eq. 1.9.31]{NIST:DLMF} to obtain an integral representation for the coefficients $c_{j}$. If $C_r$ is the positively orientated circle $\{z:|z|=r\}$ for arbitrary $r \in (0,1)$ then from (\ref{1.8.1}) and (\ref{1.9.3})
\begin{equation}
c_{j}={\frac{1}{2\pi i}\oint_{C_r}\frac{(1+z)^{1/z}}{z^{j+1}}\,dz=\frac{1}{
2\pi i}\oint_{C_r}\frac{\exp \left\{ z^{-1}\ln (1+z)\right\} }{z^{j+1}}\,dz}.
\label{1.9.7}
\end{equation}

In the second integral of (\ref{1.9.7}) we rewrite the exponential term using the geometric series 
\begin{equation}
\frac{1}{z}=-\frac{1}{1-(1+z)}
=-(1+\delta +\delta^{2}+\delta ^{3}+\cdots ), 
\label{1.9.8}
\end{equation}
where $\delta =1+z$, assuming $0<|\delta |<1$. So from (\ref{1.9.3}) 
\begin{equation}
f(z)=\exp \left\{ -\ln (\delta )\left( 1+\sum_{j=1}
^{\infty }\delta^{j}\right) \right\} 
= \frac{1}{\delta}\exp \left\{ -\ln (\delta
)\sum_{j=1}^{\infty }\delta ^{j}\right\} .  
\label{1.9.9}
\end{equation}

Note that $\delta ^{j}\ln (\delta )\rightarrow 0$ as $\delta \rightarrow 0$ for $j=1,2,3,\cdots $ by L'Hopital's rule. So using the Maclaurin expansion of the exponential function along with (\ref{1.9.9}) this function has the expansion
\begin{equation}
f(z)=\frac{1}{\delta}\left( 1-v+\frac{v^{2}}{2!}-\frac{v^{3}}{3!}+\cdots \right),
\label{1.9.11}
\end{equation}
for $0<|\delta |<1$, where $v=\ln (\delta )\sum_{j=1}^{\infty }\delta ^{j}$. From this one deduces for small $\delta $ that
\begin{equation}
v=-\delta \ln (\delta )-\delta ^{2}\ln (\delta )+{\mathcal{O}}\left\{ \delta
^{3}\ln (\delta )\right\},\,
v^{2}=\delta ^{2}\ln ^{2}(\delta )+{\mathcal{O}}\left\{ \delta ^{3}\ln
^{2}(\delta )\right\},
\label{1.9.13}
\end{equation}
and 
\begin{equation}
v^{j}={\mathcal{O}}\left\{ \delta ^{3}\ln ^{3}(\delta )\right\} \quad
(j=3,4,5,\cdots ).
\label{1.9.15}
\end{equation}

Recalling $\delta =1+z$ we consequently have from (\ref{1.9.11}) - (\ref{1.9.15})
\begin{equation}
f(z)=\left( 1+z\right) ^{-1}-\ln \left( 1+z\right) +R(z),
\label{1.9.15a}
\end{equation}
where
\begin{multline}
R(z)=(1+z)^{1/z}-\left( 1+z\right) ^{-1}+\ln (1+z) \\
=\left[ \tfrac{1}{2}\ln \left( 1+z\right) ^{2}-\ln \left( 1+z\right) \right]
\left( 1+z\right) +{\mathcal{O}}\left\{ \ln \left( 1+z\right) ^{3}\left(
1+z\right) ^{2}\right\} ,
\label{1.9.16}
\end{multline}
as $z\rightarrow -1$.

Now substitute (\ref{1.9.15a}) into (\ref{1.9.7}) to get
\begin{equation}
c_{j}=I_{1,j}+I_{2,j}+\eta_{j},
\label{I0}
\end{equation}
where 
\begin{equation}
I_{1,j}=\frac{1}{2\pi i}\oint_{C_r}\frac{1}{z^{j+1}(1+z)}\,dz,\,
I_{2,j}=-\frac{1}{2\pi i}\oint_{C_r}\frac{\ln \left( 1+z\right) }{z^{j+1}}\,dz,
\label{I1}
\end{equation}
and
\begin{equation}
\eta_{j}=\frac{1}{2\pi i}\oint_{C_r}\frac{R(z)}{z^{j+1}}\,dz.
\label{1.10.4}
\end{equation}

The integrals in (\ref{I1}) can readily be evaluated by residue theory. For the first we have by the geometric series expansion 
\begin{equation}
I_{1,j}=\underset{z=0}{\mathrm{Res}}
\left\{\frac{1}{z^{j+1}(z+1)}\right\} 
=\underset{z=0}{\mathrm{Res}}
\left\{ \sum_{s=0}^{\infty }(-1)^{s}z^{s-j-1}\right\} ={(-1)^{j}}.
\label{1.9.20}
\end{equation}
Likewise for $I_{2,j}$ one finds that
\begin{equation}
I_{2,j}=-\underset{z=0}{\mathrm{Res}}
\left\{\frac{\ln(1+z)}{z^{j+1}} \right\} = \frac{(-1)^{j}}{j}.
\label{1.9.21}
\end{equation}

\begin{figure}[ht!]
  \centering
    \includegraphics[width=130mm, scale=1, trim=0 270 0 220, clip]{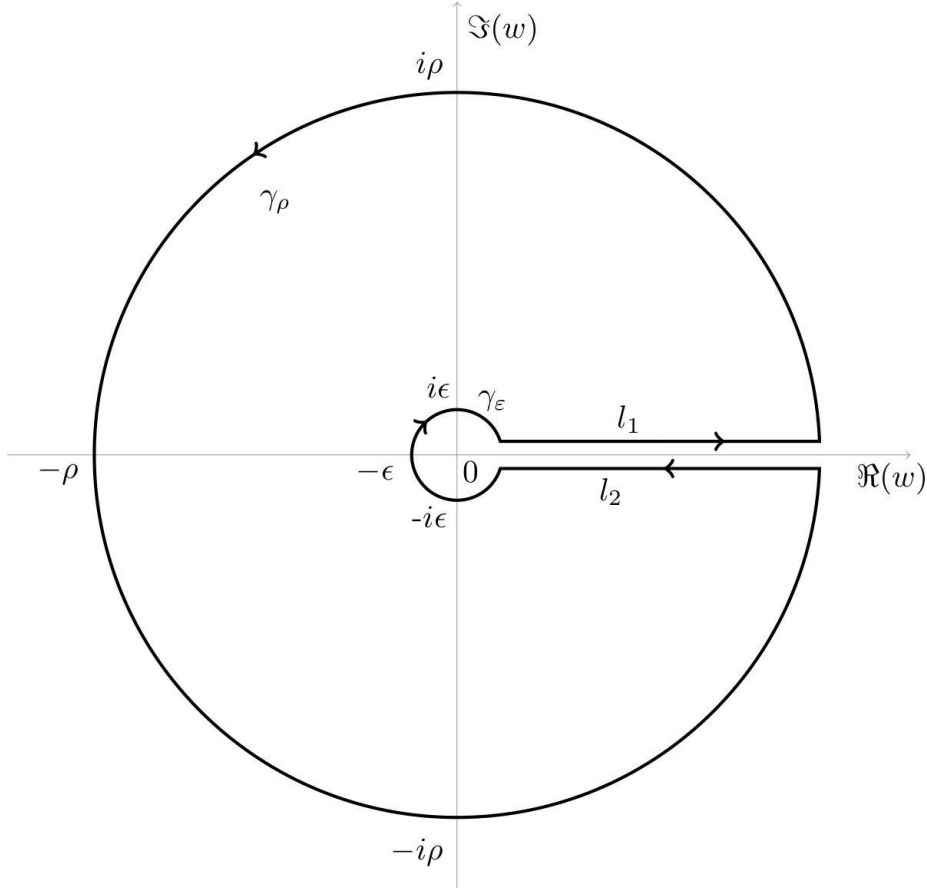}
    \caption{Contour $\Gamma_{\epsilon,\rho}$ \label{fig:fig1}}
    \label{fig1}
\end{figure}

For the integral (\ref{1.10.4}) we make a change of variable $w=-\left( z+1\right)$ to obtain the following
\begin{equation}
\eta_{j}=\frac{(-1)^{j+1}}{2\pi i}\oint_{C_r^{\prime }}
\frac{R(-1-w)}{(1+w)^{j+1}}\,dw.
\label{1.10.7}
\end{equation}
The contour $C_r^{\prime }$ in the $w$ plane is now the circle $\{w:|w+1|=r\}$ for $0<r<1$, and is positively orientated. This lies in the left half plane and encircles $w=-1$. The integrand of (\ref{1.10.7}) has a branch point at $w=0$ and a pole at $w=-1$, and is analytic elsewhere in the $w$ plane having a cut along the non-negative real axis. So we can deform the contour to a new one, called $\Gamma_{\epsilon,\rho}$, as seen in \cref{fig1}.


This contour consists of circles $\gamma _{\epsilon
}$ and $\gamma _{\rho }$\ centered at $w=0,$ of radius $\epsilon $ and $\rho $ (respectively), where $0<\epsilon <1<\rho $, and horizontal line segments $l_{1}$ and $l_{2}$ with end points $w=\epsilon \pm i0$ and $w=\rho \pm i0$ above and below the cut.

Now from (\ref{1.9.16}) we see that the integrand of (\ref{1.10.7}) is $\mathcal{O}\{w\ln(w)^{2}\}$ as $w\rightarrow 0$, and $\mathcal{O}\{w^{-j-1}\ln (w)\}$ as $w\rightarrow \infty $. Hence the contributions of $\gamma _{\epsilon }$ and $\gamma_{\rho }$ vanish as $\epsilon \rightarrow {0}$ and $\rho \rightarrow \infty$ . Then the only contribution will be along $l_{1}\cup l_{2}$, where now these lines extend from $0$ to $\infty $. We are therefore left with
\begin{equation}
\eta_{j}=\frac{(-1)^{j+1}}{2\pi i}\int_{l_{1}\cup l_{2}}\frac{R(-1-w)}{
(1+w)^{j+1}}\,dw.  
\label{1.10.8}
\end{equation}

Next, the contributions of real terms in the integrand of (\ref{1.10.8}) cancel. Hence, using $\Im\{\ln(z)^2\}=2\arg(z)\ln(|z|)$ ($z\in \mathbb{C}\setminus\{0\}$), we have from (\ref{1.9.16}) $\Im\{R(-1-w)\}=\mathcal{O}\{w\ln (w)\}$ uniformly for $w\in l_{1}\cup l_{2}$, and so from (\ref{1.10.8}) for unbounded $j$
\begin{equation}
\eta_{j}=\mathcal{O}\left\{ \int_{0}^{\infty } 
\frac{w \ln(w)}{(1+w)^{j+1}}dw\right\}.  
\label{1.10.10}
\end{equation}
We then let $1+w=e^{t}$, and consequently from (\ref{1.10.10}) obtain
\begin{equation}
\eta_{j}=\mathcal{O}\left\{ \int_{0}^{\infty }{e^{-jt}(e^{t}-1)\ln
(e^{t}-1)\,dt}\right\}.
\label{1.10.12}
\end{equation}

Now split the integral in (\ref{1.10.12}) into two integrals, one from $t=0$ to $t=1$ and the other from $t=1$ to $t=\infty$. For the first  use $e^{t}-1=\mathcal{O}(t)$ for $0\leq t\leq 1$, and for the second  use $\ln (e^{t}-1)=\mathcal{O}(t)$ for $1\leq t<\infty $. Thus as $j\rightarrow {\infty}$ we deduce that
\begin{multline}
\eta_{j}=\mathcal{O}\left\{ \int_{0}^{1}{e^{-jt}t
\ln (t)\,dt}\right\}
+\mathcal{O}\left\{ \int_{1}^{\infty }e^{-(j-1)t}t\,dt\right\}  \\
=\mathcal{O}\left( \frac{\ln (j)}{j^2}\right)
+\mathcal{O}\left(\frac{1}{je^{j}}\right)
=\mathcal{O}\left( \frac{\ln (j)}{j^2}\right),
\label{1.10.13}
\end{multline}
where the third $\mathcal{O}$ term comes from \cite[Chap. 9, Thm. 1.1]{Olver:1997:ASF}, and the fourth $\mathcal{O}$ term came from integration by parts. Finally, from (\ref{I0}), ( \ref{1.9.20}), (\ref{1.9.21}) and (\ref{1.10.13}) we arrive at (\ref{main3}).

\section*{Addendum} Christian Berg kindly made us aware of his recent joint paper \cite{Berg:2021:FEF}. In this they study a more general function $h_{\alpha}(z)=(1+1/z)^{\alpha z}$ where $\alpha>0$. They obtain a number of results for the coefficients in a Maclaurin series in $z$ for $e^{-\alpha}h_{\alpha}(-1/z)$, which are polynomials in $\alpha$. In particular they use properties of the exponential Bell partition polynomials to obtain a recursion relation, which for $\alpha=1$ reduces to our formula (\ref{main2}). In comparison we only require elementary techniques to obtain (\ref{main2}).

\section*{Acknowledgments}
TMD acknowledges financial support from Ministerio de Ciencia e Innovaci\'on, Spain, 
projects MTM2015-67142-P (MINECO/FEDER, UE) and PGC2018-098279-B-I00 (MCIU/AEI/FEDER, UE). 

\bibliographystyle{siamplain}
\bibliography{biblio}
\end{document}

%% file: ex_shared.tex
\usepackage{lipsum}
\usepackage{amsfonts}
\usepackage{graphicx}
\usepackage{epstopdf}
\usepackage{algorithmic}
\ifpdf
  \DeclareGraphicsExtensions{.eps,.pdf,.png,.jpg}
\else
  \DeclareGraphicsExtensions{.eps}
\fi

\newsiamremark{hypothesis}{Hypothesis}
\crefname{hypothesis}{Hypothesis}{Hypotheses}
\newsiamthm{claim}{Claim}
\newtheorem{remark}{Remark}
\nolinenumbers

\makeatletter
\g@addto@macro\bfseries{\boldmath}
\makeatother

\headers{coefficients in an asymptotic expansion}{T. M. Dunster and Jessica M. Perez}

\title{On the coefficients in an asymptotic expansion of $(1+1/\MakeLowercase{x})^{\MakeLowercase{x}}$}

\author{T. M. Dunster\thanks{Department of Mathematics and Statistics, San Diego State University, 5500 Campanile Drive, San Diego, CA 92182, USA. 
  (\email{mdunster@sdsu.edu}, \url{https://tmdunster.sdsu.edu}).}
\and Jessica M. Perez\thanks{Department of Mathematics and Statistics, San Diego State University, 5500 Campanile Drive, San Diego, CA 92182, USA. 
  (\email{jessica.perez.aof@gmail.com}). }
 }

\usepackage{amsopn}

\usepackage{lipsum}
\usepackage{mathrsfs}
\usepackage{amsfonts}
\usepackage{graphicx}
\usepackage{epstopdf}
\usepackage{algorithmic}
\usepackage{amssymb}
\usepackage{amsmath}

\makeatletter
\newcommand*{\addFileDependency}[1]{
  \typeout{(#1)}
  \@addtofilelist{#1}
  \IfFileExists{#1}{}{\typeout{No file #1.}}
}
\makeatother


%% file: Coefficients.bbl
\begin{thebibliography}{1}

\bibitem{Berg:2021:FEF}
{\sc C.~Berg, E.~Massa, and A.~P. Peron}, {\em A family of entire functions
  connecting the {B}essel function {$J_{1}$} and the {L}ambert {W} function},
  Constr. Approx., 53 (2021), pp.~121--154.

\bibitem{Brown:2014:CVA}
{\sc J.~W. Brown and R.~V. Churchill}, {\em Complex variables and
  applications}, Boston: McGraw-Hill Higher Education, 9~ed., 2014.

\bibitem{Chen:2014:AAF}
{\sc C.-P. Chen and J.~Choi}, {\em An asymptotic formula for $(1 + 1/x)^x$
  based on the partition function}, Amer. Math. Monthly, 44 (2014),
  pp.~338--343.

\bibitem{NIST:DLMF}
{\em {\it NIST Digital Library of Mathematical Functions}}.
\newblock http://dlmf.nist.gov/, Release 1.1.1 of 2021-03-15,
  \url{http://dlmf.nist.gov/}.
\newblock F.~W.~J. Olver, A.~B. {Olde Daalhuis}, D.~W. Lozier, B.~I. Schneider,
  R.~F. Boisvert, C.~W. Clark, B.~R. Miller, B.~V. Saunders, H.~S. Cohl, and
  M.~A. McClain, eds.

\bibitem{Hardy:1918:AFI}
{\sc G.~H. Hardy and S.~Ramanujan}, {\em Asymptotic formulae in combinatory
  analysis}, Proc. London Math. Soc., 17 (1918), pp.~75--115.

\bibitem{Olver:1997:ASF}
{\sc F.~W.~J. Olver}, {\em Asymptotics and special functions}, AKP Classics, A
  K Peters Ltd., Wellesley, MA, 1997.
\newblock Reprint of the 1974 original [Academic Press, New York].

\bibitem{Ponomarenko:2015:AFF}
{\sc V.~Ponomarenko}, {\em Asymptotic formula for $(1+1/x)^x$, revisited},
  Amer. Math. Monthly, 122 (2015), p.~587.

\end{thebibliography}
